\documentclass[10pt]{article}

\usepackage{flushend, cuted}
\usepackage{amsmath}
\usepackage{amsfonts}
\usepackage{bbm}

\usepackage{amssymb}
\usepackage{float}
\usepackage{color,amsmath,amssymb, amsfonts,amstext,amsthm}

\usepackage{epsfig, graphicx, graphics}
\usepackage{longtable}
\usepackage{cite}
\usepackage{subfigure}

\renewcommand{\phi}{\varphi}

\title{Bifurcation in mean phase portraits for stochastic dynamical systems with  multiplicative Gaussian noise \footnote{
$\quad$$^{1}$ School of Mathematics and Statistics, Zhengzhou University, 100 Kexue Road, Zhengzhou 450001, China.
$\quad$$^{2}$ Center for Mathematical Sciences, Huazhong University of Science and Technology, Wuhan 430074, China.
$\quad$$^{3}$ Department of Applied Mathematics, Illinois Institute of Technology, Chicago 60616, USA.}}

\author{Hui Wang$^{1}$, Athanasios Tsiairis$^{2}$, Jinqiao Duan$^{3}$}

\date{\today}

\begin{document}

\bibliographystyle{plain}

\maketitle

\begin{abstract}
We investigate the   bifurcation phenomena for stochastic systems with multiplicative Gaussian noise, by examining qualitative changes in mean phase portraits.  Starting from the Fokker-Planck equation for the probability density function of  solution processes, we  compute the   mean orbits  and  mean equilibrium states.  A change in the number  or stability type, when a parameter varies,  indicates a stochastic bifurcation.  Specifically,  we study stochastic bifurcation for three prototypical  dynamical systems (i.e., saddle-node, transcritical, and pitchfork systems) under multiplicative Gaussian noise,  and have  found some   interesting phenomena in contrast to the corresponding deterministic counterparts.

\end{abstract}

\section{Introduction}

Stochastic bifurcation is a phenomenon of qualitative changes in dynamical behaviors  for stochastic  systems, when a parameter varies \cite[Ch.9]{Arnold}. One of the   driving forces for studying stochastic bifurcation is that we want to know how a deterministic bifurcation differs under the influence of noise. Stochastic bifurcations have been observed in a wide range of complex systems in physical science and engineering \cite{Horsthemke2006, Deco, Ebeling,Klafter2011,Radu}.

Despite the fact that great progress has been made in the development of stochastic dynamical systems,   the study of stochastic bifurcation is   still in its infancy \cite{Crauel,Lamb}.
A stochastic bifurcation may be defined as a qualitative change, such as  the location, number and stability of equilibrium states in the evolution of a stochastic dynamical system,  as  a system parameter varies.
Usually, a bifurcation  diagram \cite{GH, Wiggins, SH, Sri} in terms of equilibrium states  vs. bifurcation parameter   is used to depict  the  qualitative changes in the phase space orbit structures for   deterministic bifurcations in low dimensional dynamical systems. A bifurcation  diagram  thus represents the qualitative change of equilibrium states (or other geometrical invariant structures)  vividly in phase portraits. However, the phase portrait for a stochastic dynamical system is a  complicated matter, due to orbits' dependency on  random samples \cite[Ch. 5]{Duan2015}.    More background in stochastic bifurcation   is reviewed  in our previous paper \cite{WH}.


Arnold \cite[Ch.9]{Arnold} considered   stochastic bifurcations by examining changes in random attractors.  Other existing works investigate stochastic bifurcation, for example, by examining the qualitative changes in  random complete quasi-solutions \cite{WangB},  in  invariant measures and their spectral stability \cite{Crauel, Lamb},  in  invariant measures with supports \cite{Xu}, in  Conley index\cite{ChenXP}, or in  logistic
map\cite{Son}.  Note that the  invariant measures or Conley index  are  actually not in the state space (where orbits or phase portraits live).  As in deterministic bifurcation, we would like to consider stochastic bifurcation in phase portraits. In our previous paper\cite{WH}, we studied the stochastic bifurcation   by examining the qualitative changes of equilibrium states in  its most probable phase portraits\cite{Duan2015,Zhuan}. The most probable orbits indicate the most likely locations of   dynamical orbits for a stochastic system.

In this present work, we  examine   stochastic bifurcation in  mean phase portraits in state space. The mean orbits indicate the expected locations of the dynamical orbits.
Sample solution orbits in phase portraits for stochastic differential equations are unintelligible objects. The phase portraits in terms of mean orbits \cite[Ch.5]{Duan2015} offer one promising option.  Thus  we propose    to study   stochastic bifurcation by examining the qualitative changes (especially the changes in the number, location and stability type for equilibrium states)  in mean phase portraits.  Specifically, we consider   bifurcation for three  prototypical   scalar    differential equations,  i.e., saddle-node, transcritical, pitchfork systems,  with multiplicative  Brownian motion.


This letter is organized as follows. In Section 2, we  define the mean phase portraits, and discuss the numerical methods for   bifurcation diagrams. In Section 3, we present bifurcation diagrams for  stochastic bifurcation under multiplicative Brownian motion in saddle-node, transcritical and pitchfork systems.  We end this letter with a brief discussion in Section 4.

\section{Methods}

Consider a scalar stochastic differential equation with multiplicative Gaussian noise
\begin{equation} \label{sde2}
  d X_t = f(r, X_t) dt + \sigma (X_t) d B_t,  \;\;  X_0= x_0.
\end{equation}
where $f$ is a given   vector field, $r$ is a real parameter,      $\sigma$ is the noise intensity, and $B_t$ is a scalar Brownian motion.

The generator for this  stochastic differential equation is
\begin{equation}
A\varphi(x)=f(r, x)\varphi'(x)  + \frac{1}{2 }\sigma^2(x) \varphi''(x) .    \label{gener1}
\end{equation}
The Fokker-Planck  equation for  the probability density function $p(x,t) \triangleq  p(x,t|x_0,0)$ of the solution process $X_t $ with initial condition $X_0=x_0$ is \cite{Duan2015}
\begin{equation} \label{fpe}
 p_t  = A^* p,  \;\;\;\;    p(x,0)=\delta(x-x_0),
\end{equation}
where  $A^*$ is the adjoint operator of the generator  $A$  in   Hilbert space $ L^2(\mathbb{R}) $, and $\delta$ is the Dirac delta function.
More specifically,
\begin{equation} \label{fpe3}
p_t =-(f(r, x)p(x, t))_x +  \frac {1}{2} (\sigma^2(x) p(x, t))_{xx},  \;\;\;   p(x,0)=\delta(x-x_0).
\end{equation}
We use a finite difference method \cite{Gao2016} to simulate this Fokker-Planck   equation (\ref{fpe3}) .
The mean orbit starting at this initial point $x_0$ in the state space $\mathbb{R}$  is then computed by
\begin{equation}\label{mean1}
 \bar{X}(x_0,t)=\int_\mathbb{R} \xi \;  p(\xi,t|x_0,0)  d\xi.
\end{equation}

A mean equilibrium state is a state which either attracts or repels  all nearby mean orbits. When it attracts    all nearby mean orbits, it is called a   mean \emph{stable} equilibrium state, while if it repels all nearby mean orbits, it is called a   mean \emph{unstable} equilibrium state. The mean phase portrait is composed of  representative mean orbits, including mean equilibrium states.


Both mean phase portraits and mean equilibrium states are  deterministic  geometric objects.  As in the study of bifurcation for deterministic dynamical systems \cite{GH,Wiggins,SH}, we examine the qualitative changes in the mean phase portraits, when a parameter varies.  A simple stochastic bifurcation is the change in the `number', `location' or `stability type' of  mean equilibrium states in the mean phase portraits.

\section{Results}

In this section, we treat explicitly  one-dimensional elementary and classical local bifurcations, namely stochastic saddle-node, transcritical  and  pitchfork bifurcation, in mean phase portraits. We  generate bifurcation diagrams by examining the mean equilibrium states for these three  prototypical stochastic systems under Gaussian noise, as  a parameter $r$ in vector field varies.

 As the rigorous general results  for mean equilibrium states  are lacking, we conduct numerical simulations to demonstrate  the  stochastic bifurcation phenomena.

\subsection{Saddle-Node Bifurcation: $d X_t = (r + X_t^2 )dt + X_t d B_t $}

First consider the stochastic saddle-node system $d X_t = f(r, X_t) dt + \sigma (X_t) d B_t$,
with $f(r, X_t)= r + X_t^2$ and $\sigma (X_t)=X_t$.
\begin{center}
\begin{figure}[th]
\subfigure[]{ \label{Fig.sub.11}
\includegraphics[width=0.45\textwidth]{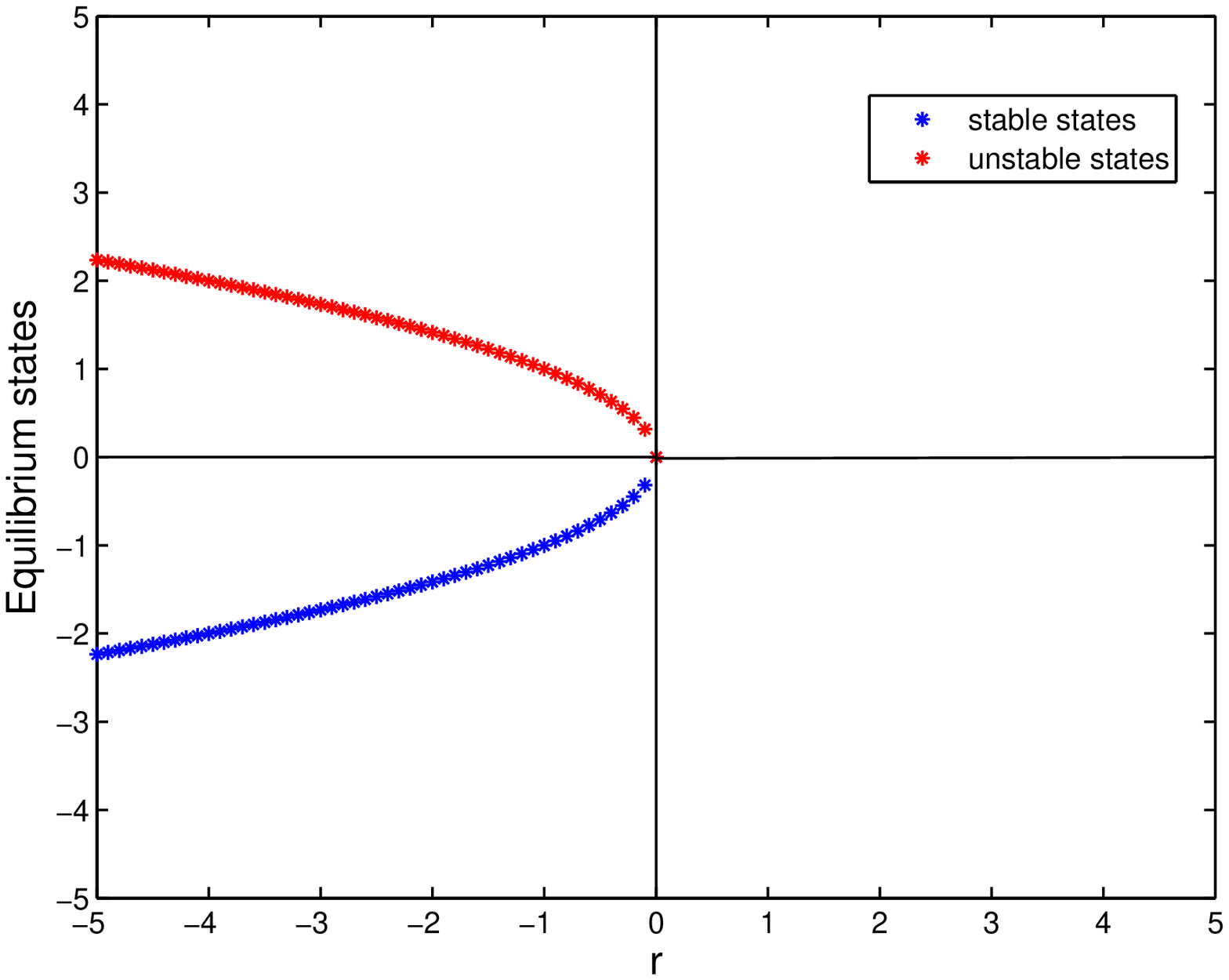}}
\subfigure[]{ \label{Fig.sub.12}
\includegraphics[width=0.45\textwidth]{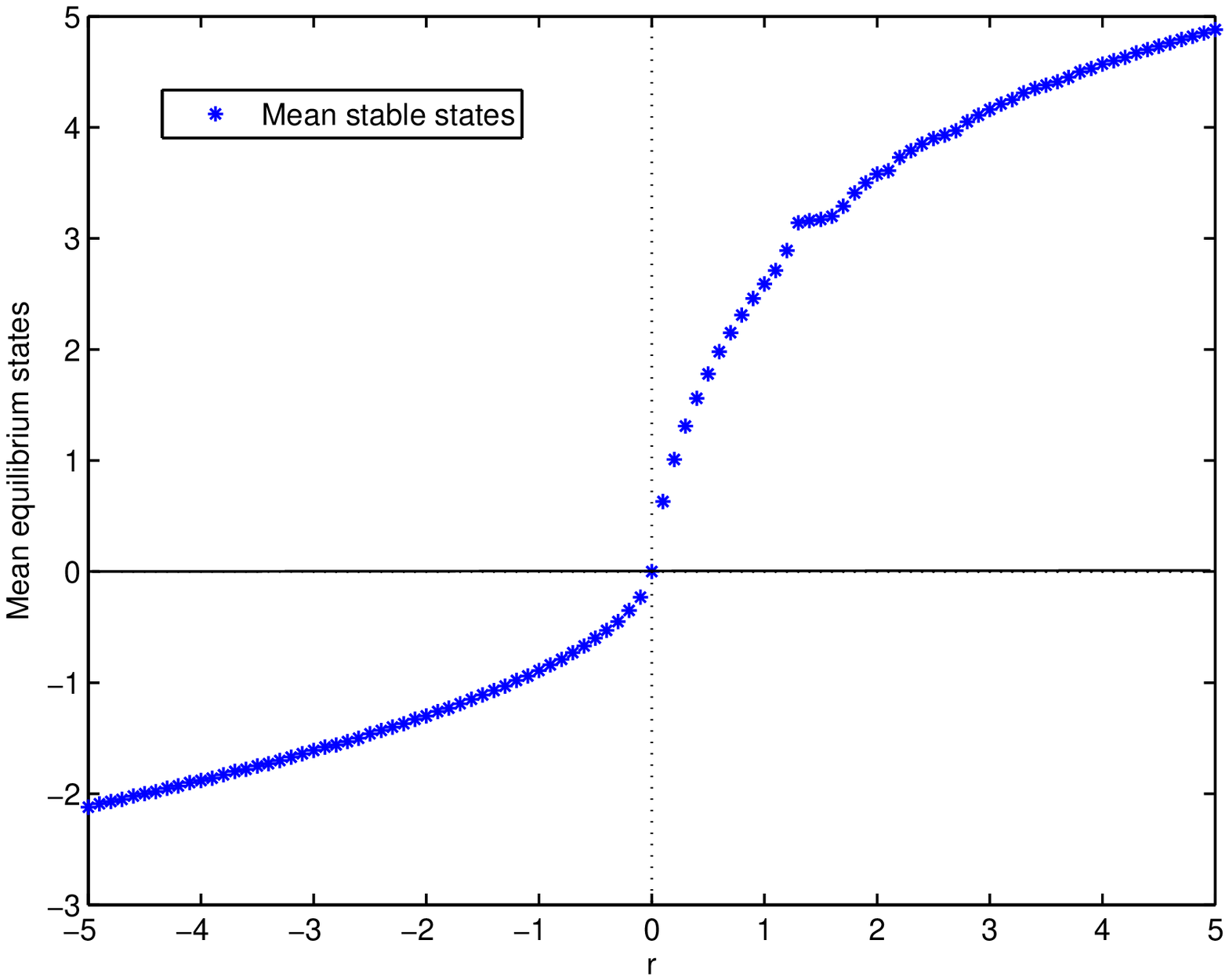}}
\caption{(Color online) Bifurcation diagrams for the saddle-node system: (a) Deterministic case (noise is absent); (b) Stochastic case. }
 \label{Fig.a1}
\end{figure}
\end{center}
Figure  \ref{Fig.sub.11} is the  well-known bifurcation diagram for the deterministic saddle-node system, i.e., equilibrium states vs. parameter $r$.  Figure  \ref{Fig.sub.12} is the bifurcation diagram for the stochastic system and it is the mean equilibrium states vs. parameter $r$.

We see that in the deterministic case, a bifurcation occurs at $r=0$. When $r$ is negative, there are two   equilibrium states, one stable and one unstable; when $r=0$ the equilibrium states become a half stable   equilibrium state; and when $r>0$, there are no   equilibrium state at all. However, in the stochastic case, for $ r>0 $, there exists one positive mean stable state;  for $ r<0 $, there exists one negative mean stable state;  while for $ r=0 $, one  mean equilibrium state is $x=0$. The location of mean equilibrium state varies as  $r$ varies. A main difference occurs at $ r>0 $, where positive mean stable states emerge.


\subsection{Transcritical Bifurcation: $d X_t = (r X_t-X_t^2) dt + X_t d B_t $}

Now consider the stochastic transcritical system $d X_t = f(r, X_t) dt + \sigma (X_t) d B_t$,
with $f(r, X_t)= r X_t- X_t^2$ and $\sigma (X_t)=X_t$.
\begin{center}
\begin{figure}[th]
\subfigure[]{ \label{Fig.sub.21}
\includegraphics[width=0.45\textwidth]{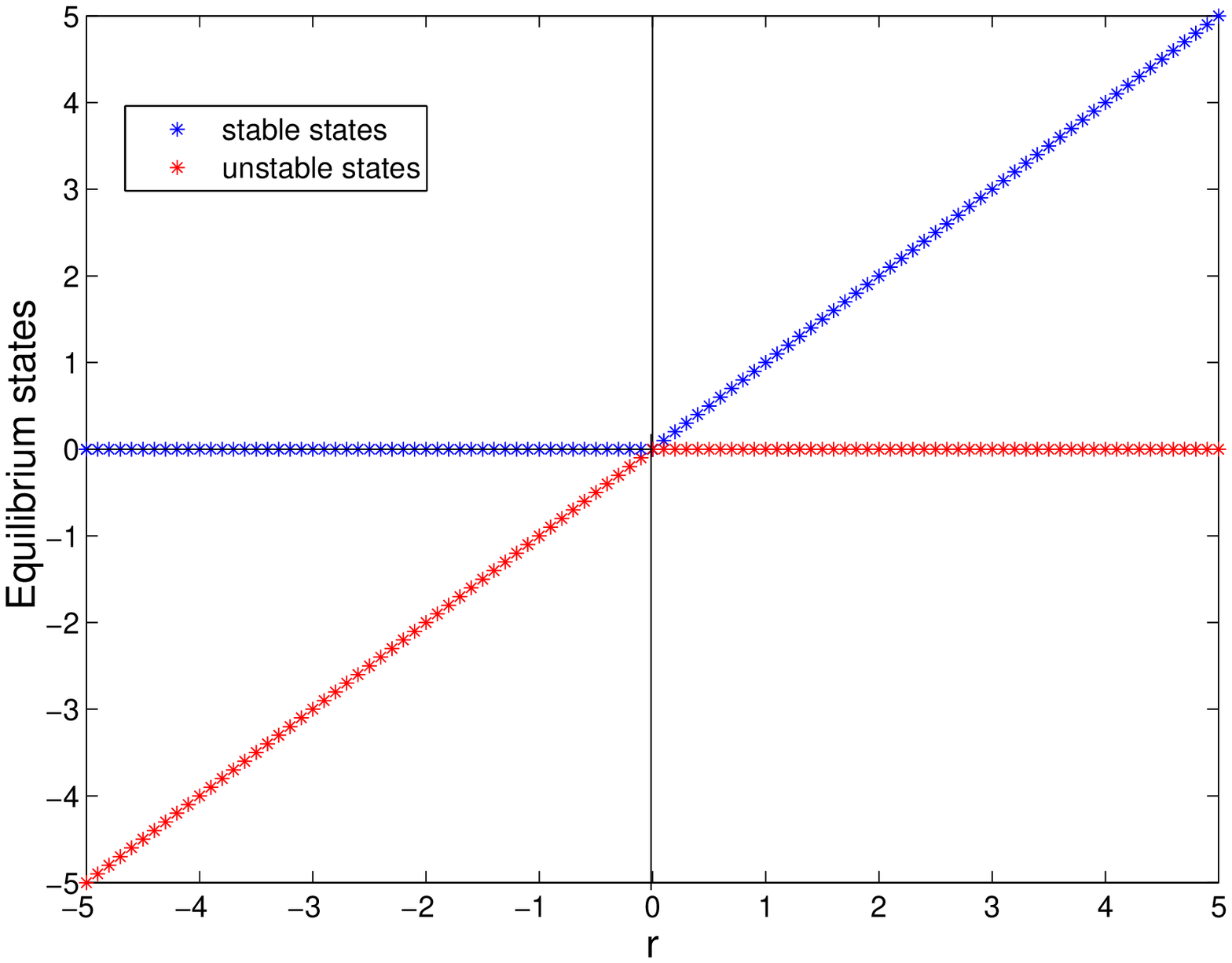}}
\subfigure[]{ \label{Fig.sub.22}
\includegraphics[width=0.45\textwidth]{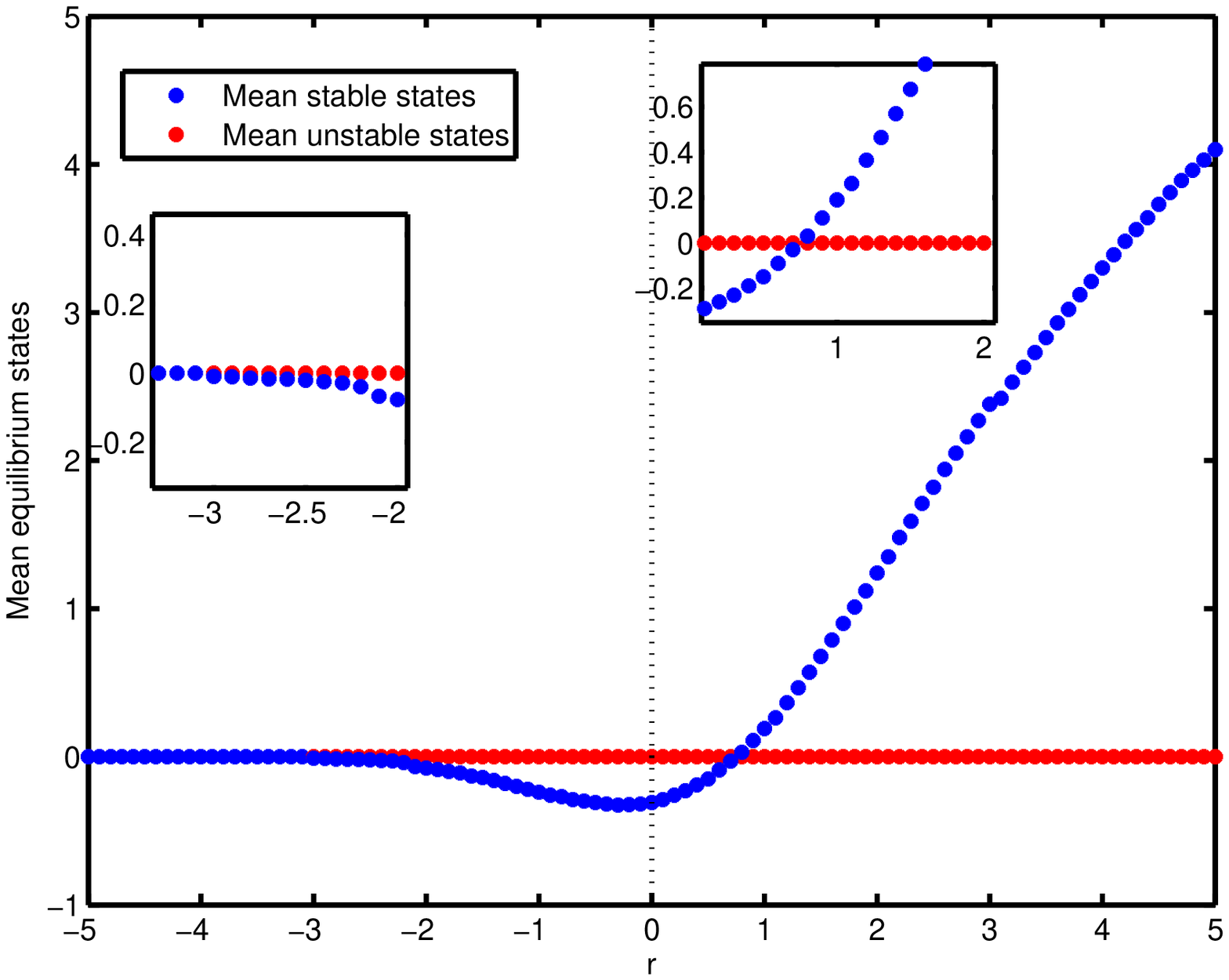}}
\caption{(Color online) Bifurcation diagrams for the transcritical  system: (a) Deterministic case (noise is absent); (b) Stochastic case.  }
 \label{Fig.a2}
\end{figure}
\end{center}

Figure  \ref{Fig.sub.21} is the bifurcation diagram for the deterministic transcritical system. Figure  \ref{Fig.sub.22} is the bifurcation diagram for the stochastic system :$d X_t = f(r, X_t) dt + \sigma (X_t) d B_t$, $f(r, X_t)= r X_t-X_t^2$, $\sigma (X_t)=X_t$.

In the deterministic case, it is the famous transcritical bifurcation.  There is an equilibrium state at $x^*=0$ for all values of $r$. For $ r<0 $, there are two equilibrium states, one unstable   at $x^*=r$ and the other  stable   at $x^*=0$;  for $r=0$ , there is only one half stable equilibrium state $x^*=0$; while for $ r>0 $, there are also two equilibrium states, one unstable   at $x^*=0$ and the other one stable   at $x^*=r$. The two equilibrium states do not disappear after the bifurcation, they just switch their stability types.

In the stochastic case, the bifurcation is very different. When $r\lesssim -3$, there  is only one mean stable point $x=0$, but for $r\gtrsim-3$, there are two mean equilibrium points: one mean stable point and one mean unstable point $0$. More specifically, the mean stable point is negative for$ -3\lesssim r \gtrsim 0.8$, while the mean stable point is positive for $r\gtrsim 0.8$. The main difference occurs at $r\lesssim-3$, where there  exists only one mean stable point $x=0$. From the details with enlarged scale in \ref{Fig.sub.22}, we can see very clearly what happens  around the bifurcation value.

\subsection{Pitchfork Bifurcation: $d X_t = (r X_t-X_t^3) dt + X_t d B_t $}
Finally, we investigate the bifurcation for the  stochastic  pitchfork system   $d X_t = f(r, X_t) dt + \sigma (X_t) d B_t$,
with $f(r, X_t)= r X_t- X_t^3$ and $\sigma (X_t)=X_t$.

\begin{center}
\begin{figure}[th]
\subfigure[]{ \label{Fig.sub.31}
\includegraphics[width=0.45\textwidth]{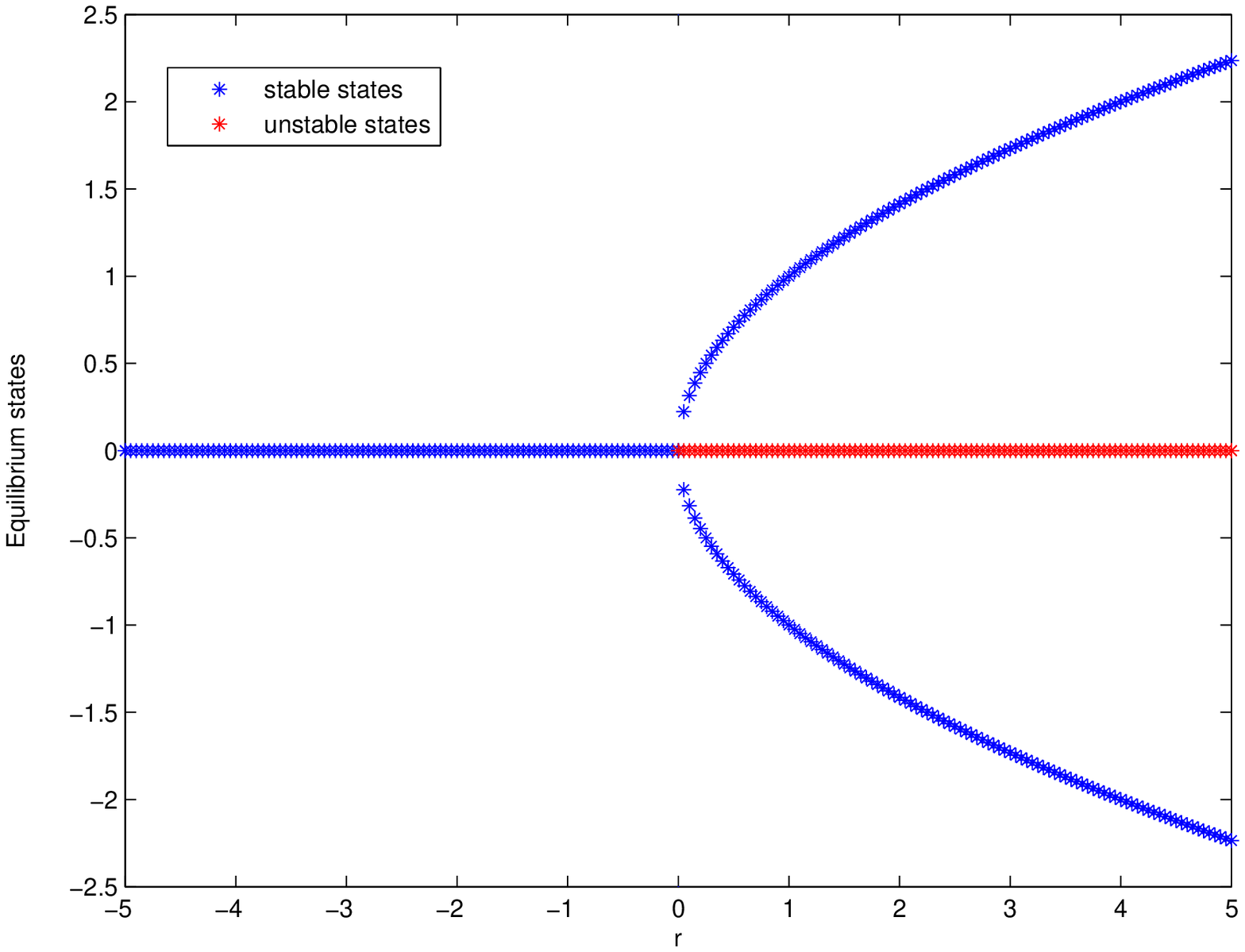}}
\subfigure[]{ \label{Fig.sub.32}
\includegraphics[width=0.45\textwidth]{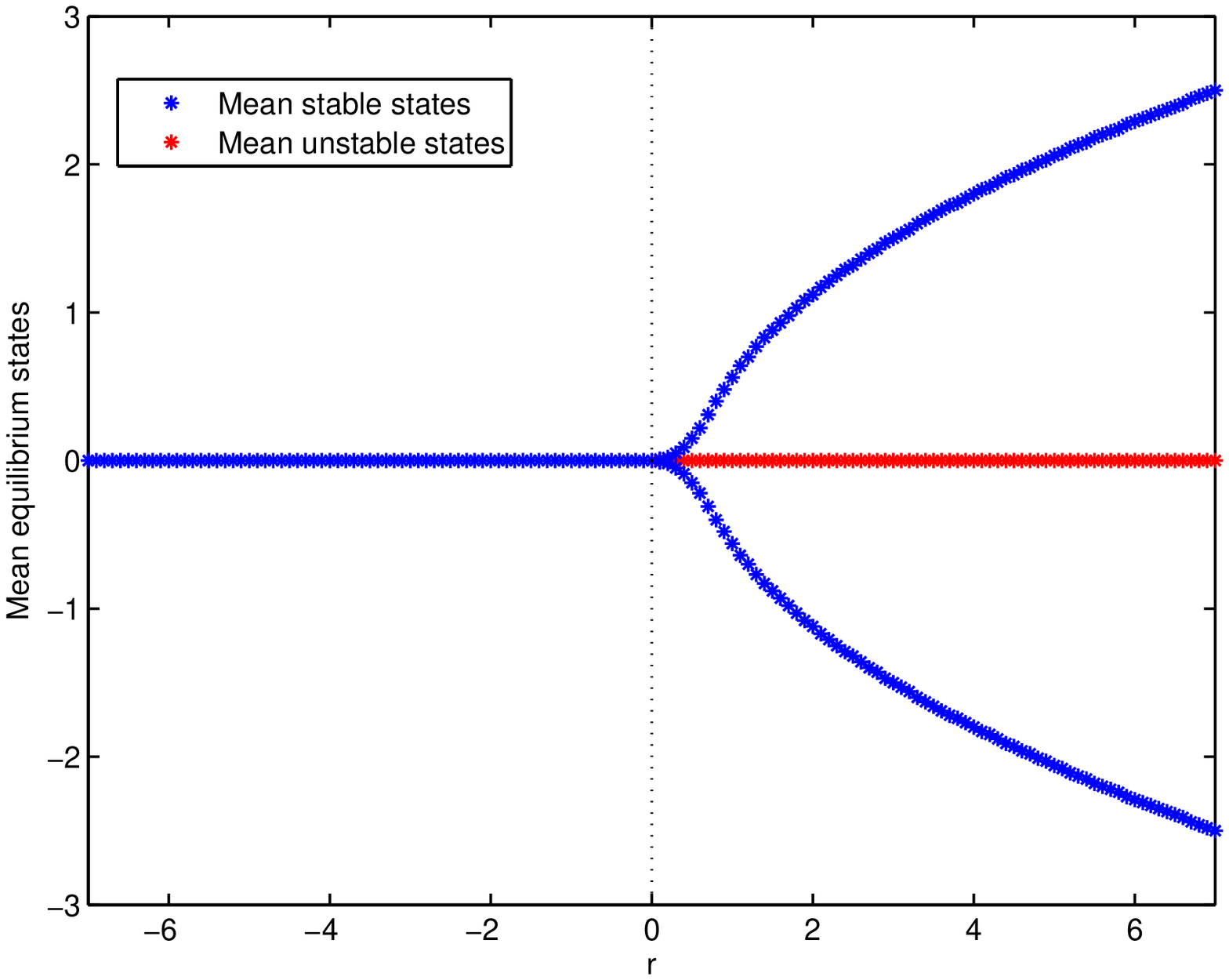}}
\caption{(Color online) Bifurcation diagrams for  the pitchfork system: (a) Deterministic case (noise is absent); (b) Stochastic case.  }
 \label{Fig.a3}
\end{figure}
\end{center}

Figure  \ref{Fig.sub.31} is the bifurcation diagram for the deterministic pitchfork system. Figure  \ref{Fig.sub.32} is the bifurcation diagram for the stochastic system:$d X_t = f(r, X_t) dt + \sigma (X_t) d B_t$, with $f(r, X_t)= r X_t-X_t^3$ and $\sigma (X_t)=X_t$.

Figure  \ref{Fig.sub.31} shows the deterministic  pitchfork system. For $r \leq 0$, $x=0$ is the only  equilibrium state which is stable. While for $ r>0 $, there exist two stable equilibrium states $ \sqrt r$ and $ -\sqrt r$ and one  unstable equilibrium state $x=0$.  The bifurcation parameter value is  $r=0$.

Figure  \ref{Fig.sub.32} is for a type of stochastic pitchfork bifurcation, with the bifurcation value     $r\approx 0.1$.  It looks     qualitatively  similar to  the deterministic one, but with a main difference in the bifurcation value,  which is ``delayed" due to the effects of noise.

\section{Conclusion}

For deterministic dynamical systems, we  usually use  their phase portraits to detect  bifurcation. But for   stochastic dynamical systems, the   phase portraits  are more delicate objects.

In this letter, we discuss stochastic bifurcation by examining the qualitative changes in mean equilibrium states
in mean phase portraits.  The mean orbits are computed by solving the corresponding Fokker-Planck equations for stochastic differential equations.

To demonstrate this stochastic bifurcation approach,  we   consider  three  prototypical     systems under multiplicative Gaussian noise. We  have observed  that   the deterministic bifurcations may be ``altered"  or  ``delayed"    by noise. It is suggested  that mean phase portraits,  composed of representative mean orbits (especially mean equilibrium states),  are   intuitive and effective tools to detect   stochastic bifurcation.

\section*{Acknowledgements}

We would like to thank Xiaoli Chen  for helpful discussions.


%

\end{document}